\documentstyle[a4,leqno,amsfonts,12pt]{article}
\newcommand{\C}{{\bf C}}

\newcommand{\OC}{\overline{{\bf C}}}

\newcommand{\N}{{\bf N}}

\newcommand{\bP}{{\bf P}}

\newcommand{\de}{\delta}
\newcommand{\La}{\Lambda}

\newcommand{\mb}{\mbox}

\newcommand{\beq}{\begin{equation}}
\newcommand{\eeq}{\end{equation}}
\newcommand{\oge}{\succeq}
\newcommand{\ole}{\preceq}
\newcommand{\ve}{\varepsilon}

\newcommand{\ov}{\overline}
\newcommand{\al}{\alpha}
\newcommand{\be}{\beta}
\newcommand{\Om}{\Omega}
\newcommand{\om}{\omega}
\newcommand{\z}{\zeta}

\newcommand{\kap}{\mb{ cap}}
\newcommand{\ga}{\gamma}

\newtheorem{th}{Theorem}

\newtheorem{lem}{Lemma}

\newcommand{\ueberschrift}{\bigskip\goodbreak\noindent\bigskip}
\newcounter{theabsatz}
\newcommand{\absatz}[1]{\stepcounter{theabsatz} \ueberschrift
                           {\large \bf \arabic{theabsatz}. {#1}} \setcounter{equation}{0}}

\begin{document}
\mathsurround=2pt

\begin{center}
{\bf \large On Hilbert lemniscate theorem
for a system of continua}\\[2ex]

V. V. ANDRIEVSKII \\[2ex]

 {\bf Abstract}

\end{center}

Let $K$ be
 a compact set in the complex plane consisting of a finite
number of continua. We study the rate of approximation of $K$ from
the outside by lemniscates in terms of level lines of the Green function
for the complement of $K$.



\absatz{Introduction and main results}

Let $K\subset\C$ be a compact set in the complex plane $\C$
consisting of disjoint  closed connected sets (continua) $K^j,j=1,2,\ldots,\nu$,
i.e.,
$$
K=\bigcup_{j=1}^\nu K^j;\quad K^j\bigcap K^k=\emptyset\mb{ for }j\neq k;
\quad\mb{ diam}(K^j)>0;
$$
where
$
\mb{diam}(S)$ 
is the diameter of $S\subset\C$.  
We always assume that $\Om:=\OC\setminus K$ is connected.
Here, $\OC:=\C\cup\{\infty\}$ is the extended complex plane.

According to the Hilbert lemniscate theorem (see \cite[p. 159]{ran}), for any open 
neighborhood $U$ of $K$, there exists a polynomial $p$ such that
\beq\label{1.1}
\frac{|p(z)|}{||p||_K}>1,\quad z\in \C\setminus U,
\eeq
where $||f||_S$ denotes the uniform norm of $f:S\to\C$ on $S\subset\C$.
Certainly, the degree of $p$ depends on $U$. 

Let
$\bP_n, n\in \N:=\{1,2,\cdots\}$ be the set of all polynomials of
degree at most $n$. Denote by $g_\Om(z)$
the Green function for $\Om$ with pole at $\infty$.
It will be convenient for us to extend the Green function to $K$ by setting it
equal to zero there.
Let $s_n(K),n\in\N$ be the infimum of $s>0$ for which there exists $p\in\bP_n$
such that (\ref{1.1}) holds with
$$
U=U_s:=\{z:g_\Om(z)<s\}.
$$
A result by Siciak \cite[Theorem 1]{sic} for the Fekete polynomials
yields that
\beq\label{1.2}
s_n(K)=O\left(\frac{\log n}{n}\right)\quad \mb{as }n\to\infty
\eeq
(cf. \cite[Theorem 1]{and00}, \cite[Theorem 2]{kit}, \cite[Theorem 2.2]{pri}).

\begin{th}\label{th1}
Under the above assumptions,
\beq\label{1.1r}
s_n(K)=O\left(\frac{(\log\log n)^2}{n}\right)\quad \mb{as }n\to\infty.
\eeq
\end{th}
We also would like to demonstrate that
if more information is known about the geometry of $K$, (\ref{1.1r})
can be further improved in the following
way.
A Jordan curve  $L\subset\C$ is called {\it quasiconformal}
(see \cite{ahl}, \cite[p. 100]{lehvir} or \cite{geh}) if
for every $z_1,z_2\in L$,
$$
\mb{diam}(L(z_1,z_2))\le \La_L|z_2-z_1|,
$$
where $L(z_1,z_2)$ is the smaller subarc of $L$ between $z_1$ and $z_2$;
 and $\La_L\ge1$ is a constant that depends only on $L$.
A {\it quasidisk} is a Jordan domain bounded by the quasiconformal curve.
\begin{th}\label{th2} If each  $K^j $ is a closed quasidisk, then
\beq\label{1.4}
s_n(K)=O\left(\frac{1}{n}\right)\quad \mb{as }n\to\infty.
\eeq
\end{th}
See \cite[Theorem 2]{and00} for a special case of 
this result.

Our proof for Theorem \ref{th2}
yields insights that can be leveraged to obtain other results.
Specifically,
 that
for sufficiently large $n$ there exists polynomial $p_n\in\bP_n$
such that  $G(p_n):=\{z:|p_n(z)|\le1\}$ consists of exactly $\nu$
disjoint Jordan domains and
$$
K\subset G(p_n)\subset U_{C/n}
$$
holds with a constant $C=C(K)>0$. Moreover, using reasoning from
\cite[Section 3]{andnaz} it can be shown that $p_n$ may be chosen such
that all its zeros belong to $K$.

It is worth pointing out that (\ref{1.4}) is optimal in the following sense. Let $K$ be a 
closed quasidisk, i.e. $\nu=1$, for which there exist $\z\in\partial K, \de>0$ and $1<\be<2$
such that a circular sector with center at $\z$, radius $\de$ and opening
$\be\pi$ is a subset of $\ov{\Om}$.
Then, according to \cite[Theorem 3]{and00},
$$
s_n(K)\ge\frac{\ve}{n},\quad n\in\N
$$
is true with some constant $\ve=\ve(K)>0$.

Furthermore, in the case $\nu=1$ it is natural to approximate $K$ by
lemniscates given by Faber polynomials $F_n=F_n(K)$.
It was shown in \cite[pp. 300-301]{and00} that for the quasidisk $K$
constructed by Gaier \cite{gai}, the inequality (\ref{1.1}) does not
hold for $p=F_n, U=U_{\al\log n/n}, $ some constant $\al=\al(K)$ and 
an infinite number of $n\in\N$ (cf. (\ref{1.2})-(\ref{1.4})).

In what follows, we use following notation. 
$$
d(S_1,S_2):=\inf_{z_1\in S_1,z_2\in S_2}|z_2-z_1|,\quad S_1,S_2\subset\C.
$$
For a Jordan curve $L\subset\C$,  denote by int$(L)$ the bounded 
connected component of $\OC\setminus L$.

For a (Borel) set $S\subset\C$, denote by $|S|$ its linear measure (length)
and by $\sigma(S)$ its two-dimensional Lebesgue measure (area).

In what follows, we denote by $c,c_1,\ldots$ positive constants
that are either absolute or they depend only on $K$. For the
nonnegative functions $a$ and $b$ we write $a\ole b$ if $a\le c_1 b$,
and $a\asymp b$ if $a\ole b$ and $b\ole a$ simultaneously.

\absatz{Construction of auxiliary polynomials}

In this section we review (in more general setting)
 the construction of  the monic polynomials suggested
in \cite{tot12, tot13, and151, andnaz}. 
For the convenience of the reader, we repeat the relevant
material from these papers without proofs, thus making our 
exposition self-contained.

We start with some general
facts from potential theory which can be found, for example,
in \cite{wid, ran, saftot}.
 The  Green function $g(z)=g_\Om(z)$  has a multiple-valued harmonic conjugate
 $\tilde{g}(z)$. Let
 $$
 \Phi(z):=\exp(g(z)+i\tilde{g}(z)),
 $$
 $$
 K_s:=\{z:g(z)=s\},\quad s>0.
 $$
 Note that
 \beq\label{2.1}
 \kap(K_s)=e^s\kap(K).
 \eeq
Here $\kap(S)$ is the logarithmic capacity of a compact set $S\subset\C$.

 Let $s_0>0$  be such that for $0<s<s_0$, the set
 $K_s=\cup_{j=1}^\nu K_s^j$ consists  of $\nu$ mutually disjoined
 Jordan curves, where $K_s^j$ is the curve surrounding $K^j$.
 Moreover, we  fix a positive number  $s^*<s_0/10$ so small  that
for each $j=1,\ldots,\nu$,
\beq\label{2.1n}
 d(\z,K^j)\le d(\z,K^j_{s_0}),\quad \z\in \mb{int}(K^j_{s^*}).
 \eeq
 Let $\mu=\mu_K$ be the equilibrium measure of $K$
   and let $\om_j:=\mu(K^j).$
The function
 $
 \phi_j:=\Phi^{1/\om_j}(\z)
 $
is a conformal and univalent mapping of
$
\Om^j:=\mb{int}(K_{s_0}^j)\setminus K^j
$
 onto the annulus
 $
 A^j:=\{w:1<|w|<e^{s_0/\om_j}\}
 $
 as well as
 $$
 K_s^j=\{\z\in \Om^j:|\phi_j(\z)|=e^{s/\om_j}\},\quad 0<s<s_0.
 $$

 Note that for $\mu_s:=\mu_{K_s} $,
 $$
 \mu_s(K_s^j)=\mu(K^j)=\om_j,\quad 0<s<s_0.
 $$
 Furthermore,  for an arc
 $$
 \ga=\{\z\in K_s^j:\theta_1\le\arg \phi_j(\z)\le\theta_2\},
 \quad 0<\theta_2-\theta_1\le 2\pi,
 $$
 we have
 $$
 \mu_s(\ga)=\frac{(\theta_2-\theta_1)\om_j}{2\pi}\, .
 $$
 Assuming that $m\in\N$ is sufficiently large, i.e.
 $m>10/(\min_{j}\om_j)$ we let
 $$
 m_j:=\lfloor m \om_j\rfloor,\quad j=1,2,\ldots,\nu-1,
 $$
 $$
 m_\nu:=m-(m_1+\ldots+m_{\nu-1}),
 $$
 where $\lfloor a\rfloor$ means the integer part of  a real number $a$.

 Therefore,
 \beq\label{2.6}
 0\le m_\nu-m\om_\nu=\sum_{j=1}^{\nu-1}(m\om_j-m_j)\le \nu-1.
 \eeq
 Next, for $0<s<s^*$, we represent each $K^j_s$
 as the union of closed subarcs $I^j_{s,k},k=1,\ldots,m_j$ such that
 $$
 I_{s,k}^j\bigcap I_{s,k+1}^j=:\xi_{s,k}^j,\quad k=1,\ldots, m_j-1,
 $$
 and $I_{s,m_j}^j\cap I_{s,1}^j=:\xi_{s,m_j}^j=:\xi_{s,0}^j$ are points of $K_s^j$
 ordered in a positive  direction, as well as
 $$
 \mu_s(I_{s,k}^j)=\frac{\om_j}{m_j},\quad k=1,\ldots,m_j.
 $$
Consider also $\psi_j:=\phi_j^{-1}$,
$$
\tilde{D}^j_{s,k}:=\left\{t=re^{i\eta}:\psi_j(e^{s/\om_j+i\eta})\in I^j_{s,k},
0\le e^{s/\om_j}-r\le\frac{e^{s/\om_j}-1}{64}\right\},
$$
$$
D^j_{s,k}:=\psi_j(\tilde{D}^j_{s,k}),\quad D^j_s:=\bigcup_{k=1}^{m_j}D^j_{s,k},
\quad D_s:=\bigcup_{j=1}^\nu D_s^j.
$$
\begin{lem}\label{lem2.1}
Let $m,q\in\N$ and $c:=640\pi \max_je^{s^*/\om_j}.$
 Then
 for $m\ge m_0:=\lfloor 2cq/s^*+10\nu/\min_j\om_j\rfloor, s=cq/m<s^*, j=1,\ldots,\nu$
 and $k=1,\ldots,m_j$, we have
\beq\label{2.5h}
\frac{d(\xi^j_{s,k},K^j)}{c^2q}\le |\xi_{s,k}^j-\xi^j_{s,k-1}|\le \mb{\em diam}(I^j_{s,k})\le
|I^j_{s,k}|\le \frac{d(I^j_{s,k},K^j)}{10q}  .
\eeq
Moreover, if $q=1$ then
\beq\label{2.1du}
\sigma(D^j_{s,k})\ge\frac{d(\xi^j_{s,k},K^j)^2}{c^2}
\eeq
as well as
\beq\label{2.2du}
\mb{\em diam}(D^j_{s,k})\le \frac{1}{2}d(\xi^j_{s,k},K^j).
\eeq
\end{lem}
For the proof of Lemma \ref{lem2.1}, see Section 3.

Next, we construct the  points $\z_{s,k,l}^j,  l =1,\ldots, q$ as follows. 
For  $s=cq/m$ as in Lemma \ref{lem2.1} and $u=1,\ldots, q$, let
$$
 m_{s,k,u}^j:=\frac{1}{\mu_s(I^j_{s,k})}
\int_{I^j_{s,k}}(\xi-\xi^j_{s,k})^ud\mu_s(\xi).
$$
Consider the system of equations
$$
\sum_{l=1}^q(r_{s,k,l}^j)^u=qm_{s,k,u}^j=:\tilde{m}_{s,k,u}^j, \quad u=1,\ldots,q.
$$
We interpret $r_l:=r_{s,k,l}^j$ as the roots of the polynomial $z^q+a_{q-1}z^{q-1}+\ldots+a_0$
whose coefficients satisfy Newton's identities
\beq\label{2.1h}
\tilde{m}_u+a_{q-1}\tilde{m}_{u-1}+\ldots+a_{q-u+1}\tilde{m}_{1}=-ua_{q-u},\quad u=1,\ldots,q,
\eeq
where $\tilde{m}_u:=\tilde{m}_{s,k,u}^j$ satisfy $|\tilde{m}_u|\le q d^u,
d=d_{s,k}^j:=$diam$(I^j_{s,k})$.

According to (\ref{2.1h}),
$$
|a_{q-u}|\le q^ud^u,\quad u=1,\ldots,q,
$$
which implies
\beq\label{2.10h}
|r_{s,k,l}^j|\le 2qd_{s,k}^j.
\eeq
See \cite[Section 2]{andnaz} for more details.

Let $\z^j_{s,k,l}:=\xi^j_{s,k}+r_{s,k,l}^j.$
By virtue of (\ref{2.5h}) and (\ref{2.10h}),
\beq\label{2.2pu}
\frac{|\z_{s,k,l}^j-\xi^j_{s,k}|}{d(\xi^j_{s,k},K^j)}\le 
\frac{|r_{s,k,l}^j|}{d(I^j_{s,k},K^j)}\le \frac{1}{5},
\eeq
as well as for $\xi\in I^j_{s,k}$,
\beq\label{2.3pu}
\frac{|\xi-\xi^j_{s,k}|}{d(\xi^j_{s,k},K^j)}\le
\frac{d^j_{s,k}}{d(I^j_{s,k},K^j)}\le\frac{1}{10}.
\eeq
By \cite[p. 23, Lemma 2.3]{andbla}, which is an immediate consequence of
Koebe's one-quarter theorem, and (\ref{2.1n}), we have the following lemma.
 Recall that $\psi_j$ is defined in
$
 A^j:=\{\tau:1<|\tau|<e^{s_0/\om_j}\}
 $. Let
 $E^j:=\{w:1<|w|<e^{s^*/\om_j}\}$.
\begin{lem}\label{lem3.1h}
Let $w\in E^j,\tau\in A^j$ and $\xi=\psi_j(w),\z=\psi_j(\tau)$.
Then
\beq\label{2.1pu}
\frac{1}{4}\frac{d(\xi,K^j)}{|w|-1}\le|\psi_j'(w)|\le
4\frac{d(\xi,K^j)}{|w|-1}.
\eeq
Moreover, if either $|\tau-w|\le(|w|-1)/2$  or$|\z-\xi|\le d(\xi,K^j)/2$,
then
\beq\label{3.2}
\frac{1}{16}\frac{|\tau-w|}{|w|-1}\le\frac{|\z-\xi|}{d(\xi,K^j)}\le
16\frac{|\tau-w|}{|w|-1}.
\eeq
\end{lem}
Next, we claim that for $z\in K^j_{10s},s < s^*$,
\beq\label{2.4pu}
|z-\xi^j_{s,k}|\ge\frac{1}{2}d(\xi^j_{s,k},K^j).
\eeq
Indeed, if we assume, contrary to (\ref{2.4pu}), that
$$
|z-\xi^j_{s,k}|<\frac{1}{2}d(\xi^j_{s,k},K^j),
$$
then, according to the left-hand side of (\ref{3.2}),
\begin{eqnarray*}
e^{10s/\om_j}-e^{s/\om_j}&\le&
|\phi_j(z)-\phi_j(\xi^j_{s,k})|\\
&\le&
16(e^{s/\om_j}-1)\frac{|z-\xi^j_{s,k}|}{d(\xi_{s,k}^j,K^j)}<
8(e^{s/\om_j}-1),
\end{eqnarray*}
which contradicts to the obvious inequality
$$
e^{10x}-e^x\ge8(e^x-1),\quad x\ge 0.
$$
Hence, (\ref{2.4pu}) is proven.

A major component of the proof of (\ref{1.1r}) and (\ref{1.4}) is the polynomial
$$
P_n(z):=\prod_{j=1}^\nu\prod_{k=1}^{m_j}\prod_{l=1}^{q}(z-\z^j_{s,k,l}),\quad n=qm,
s=\frac{cq}{m}\le s^*.
$$
For $z\in K_{10s}$,
\begin{eqnarray*}
&&
 mg_{\Om_s}(z)+m\log\kap(K_s)=m\int_{K_s}\log|z-\xi|
d\mu_s(\xi)\\
&=&
\sum_{j=1}^\nu\sum_{k=1}^{m_j}\left(m-\frac{m_j}{\om_j}\right)
\int_{I^j_{s,k}}\log|z-\xi|d\mu_s(\xi)\\
&&+
\sum_{j=1}^\nu\sum_{k=1}^{m_j}\frac{1}{\mu_s(I^j_{s,k})}
\int_{I^j_{s,k}}\log|z-\xi|d\mu_s(\xi)\\
&=:&
\Sigma_1(z)+\Sigma_2(z).
\end{eqnarray*}
Since 
\begin{eqnarray*}
\int_{K_s}|\log|z-\xi||d\mu_s(\xi)&\le&
|\log\mb{diam}(K_{s_0})|+\int_{K_s} \log\frac{\mb{diam}(K_{s_0})}{|z-\xi|}d\mu_s(\xi)\\
&\le& 2|\log\mb{diam}(K_{s_0})| -g_{\Om_s}(z)-\log\kap(K_s)\ole 1,
\end{eqnarray*}
according to (\ref{2.6}), for $z\in K_{10s}$, we obtain
$$
|\Sigma_1(z)|\ole \int_{K_s}|\log|z-\xi||d\mu_s(\xi)\ole 1.
$$
Next, for the same $z\in K_{10s}$,
\begin{eqnarray*}
&& 
\log |P_n(z)|-q\Sigma_2(z)\\
&=&
\sum_{j=1}^\nu\sum_{k=1}^{m_j}\sum_{l=1}^q\left(\log|z-\z^j_{s,k,l}|
-\frac{1}{\mu_s(I^j_{s,k})}\int_{I^j_{s,k}}\log|z-\xi|d\mu_s(\xi)\right)\\
&=&\sum_{j=1}^\nu \frac{m_j}{\om_j}\sum_{k=1}^{m_j}\sum_{l=1}^q\int_{I^j_{s,k}}
\log\left|\frac{z-\z^j_{s,k,l}}{z-\xi}\right| d\mu_s(\xi).
\end{eqnarray*}
Moreover, (\ref{2.5h}), (\ref{2.10h}), 
(\ref{2.4pu}) and Taylor's theorem 
\cite[pp. 125-126]{ahl1} imply
\begin{eqnarray*}
\log\left(\frac{z-\z^j_{s,k,l}}{z-\xi}\right)&=&
\log\left(1-\frac{\z^j_{s,k,l}-\xi^j_{s,k}}{z-\xi^j_{s,k}}\right)
-\log\left(1-\frac{\xi-\xi^j_{s,k}}{z-\xi^j_{s,k}}\right)\\
&=&
\sum_{u=1}^q\frac{1}{u}\left(\left(\frac{\xi-\xi^j_{s,k}}{z-\xi^j_{s,k}}\right)^u
-\left(\frac{\z^j_{s,k,l}-\xi^j_{s,k}}{z-\xi^j_{s,k}}\right)^u\right)+B_{s,k,l}^j(z),
\end{eqnarray*}
where 
$$
|B_{s,k,l}^j(z)|\ole \left(\frac{2qd_{s,k}^j}{d(z,I^j_{s,k})}\right)^{q+1}.
$$
Since 
\begin{eqnarray*}
&&\sum_{l=1}^q\int_{I^j_{s,k}}\left((\xi-\xi^j_{s,k})^u-
(\z^j_{s,k,l}-\xi^j_{s,k})^u\right)d\mu_s(\xi)\\
&=&\sum_{l=1}^q\left(\frac{\om_j}{m_j}m^j_{s,k,u}-
(r^j_{s,k,l})^u\frac{\om_j}{m_j}\right)=
\frac{\om_j}{m_j}
\sum_{l=1}^q\left(m^j_{s,k,u}-
(r^j_{s,k,l})^u\right)\\
&=&\frac{\om_j}{m_j}\left(\tilde{m}^j_{s,k,u}-\sum_{l=1}^q(r^j_{s,k,l})^u\right)=0,
\end{eqnarray*}
for $z\in K_{10s}$, we have
\begin{eqnarray*}
&&\left|\log|P_n(z)|-ng_{\Om_s}(z)-n\log\kap(K_s)\right|
\le q|\Sigma_1(z)|+
|\log|P_n(z)|-q\Sigma_2(z)|\\
&\ole& q+\left|\sum_{j=1}^\nu \frac{m_j}{\om_j}\sum_{k=1}^{m_j}\sum_{l=1}^q
\int_{I^j_{s,k}}\Re\left(\sum_{u=1}^q\frac{1}{u}\left(
\left(\frac{\xi-\xi^j_{s,k}}{z-\xi^j_{s,k}}\right)^u
-\left(\frac{\z^j_{s,k,l}-\xi^j_{s,k}}{z-\xi^j_{s,k}}\right)^u\right)
\right.\right.\\
&&\left.\left.
+B_{s,k,l}^j(z)
\right)d\mu_s(\xi)\right|\\
&\le& q+\left|\sum_{j=1}^\nu \frac{m_j}{\om_j}\sum_{k=1}^{m_j}\sum_{u=1}^q
\left(\frac{1}{u(z-\xi^j_{s,k})^u}\sum_{l=1}^q
\int_{I^j_{s,k}}\left(
(\xi-\xi^j_{s,k})^u
-(\z^j_{s,k,l}-\xi^j_{s,k})^u\right)d\mu_s(\z)\right.\right.\\
&&\left.\left.
 +\sum_{l=1}^q\int_{I^j_{s,k}}B_{s,k,l}^j(z)
d\mu_s(\xi)\right)\right|
\\
&\ole& q\left(1+\sum_{j=1}^\nu \sum_{k=1}^{m_j}
\left(\frac{2qd_{s,k}^j}{d(z,I^j_{s,k})}\right)^{q+1}\right).
\end{eqnarray*}
To summarize,
according to (\ref{2.1}) and the identity
$$
g_{\Om_s}=g_\Om(z)-s,\quad z\in\C\setminus U_s,
$$
for $z\in K_{10s}, s< s^*$,
 we have
\begin{eqnarray}
&&\left|\log|P_n(z)|-ng_{\Om}(z)-n\log\kap(K)\right|\nonumber\\
\label{2.11}
&\ole&
q^2+ q\sum_{j=1}^\nu\sum_{k=1}^{m_j}
\left(\frac{2qd_{s,k}^j}{d(z,I^j_{s,k})}\right)^{q+1}.
\end{eqnarray}

\absatz{Distortion properties of $\phi_j$ and $\psi_j$}

{\bf Proof of Lemma \ref{lem2.1}.}
Only the first and the last inequalities in (\ref{2.5h})
are not trivial.
According to our assumption
$$
\frac{1}{2}\le \frac{m\om_j}{m_j}\le 2,\quad j=1,\ldots,\nu.
$$
Let  $\eta_0^j<\eta_1^j<\ldots<\eta_{m_j}^j=\eta_0^j+2\pi$
  be determined by $t^j_{s,k}:=\phi_j(\xi_{s,k}^j)=\exp(s/\om_j+i\eta_k^j)$, i.e.,
  $\eta_k^j-\eta_{k-1}^j=2\pi/m_j$.

Let $c_2:=\max_je^{s^*/\om_j}$. For $t\in\tilde{I}^j_{s,k}:=
\phi_j(I^j_{s,k})$,
\begin{eqnarray*}
|t-t^j_{s,k}|&\le& e^{s/\om_j}|e^{i\eta_k^j}-e^{i\eta_{k-1}^j}|\\
&\le&c_2\frac{4\pi}{m_j}\le\frac{1}{32}\frac{cq}{m\om_j}=\frac{1}{32}\frac{s}{\om_j}
<\frac{1}{32}\left( e^{s/\om_j}-1\right).
\end{eqnarray*}
Since (\ref{3.2}) implies for $\xi\in I^j_{s,k}$,
\beq\label{3.1pi}
\frac{|\xi-\xi^j_{s,k}|}{d(\xi^j_{s,k},K^j)}\le 16\frac{|t-t^j_{s,k}|}{e^{s/\om_j}-1}
<\frac{1}{2},
\eeq
by (\ref{2.1pu}), for $t\in\tilde{I}^j_{s,k}$, we have 
$$
|\psi'_j(t)|\le 8\frac{d(\xi^j_{s,k},K^j)}{e^{s/\om_j}-1}
\le 16\frac{d(I^j_{s,k},K^j)}{e^{s/\om_j}-1}.
$$
Hence,
\begin{eqnarray*}
|I^j_{s,k}|&\le&
e^{s/\om_j}\int_{\eta^j_{k-1}}^{\eta^j_k}|\psi'_j(e^{s/\om_j+i\theta})|d\theta
\le 16c_2\frac{d(I^j_{s,k},K^j)}{e^{s/\om_j}-1}(\eta^j_k-\eta^j_{k-1})\\
&\le&
\frac{32\pi c_2}{cq}\frac{\om_jm}{m_j}d(I^j_{s,k},K^j)
\le	\frac{64\pi c_2}{cq}d(I^j_{s,k},K^j)\le
\frac{d(I^j_{s,k},K^j)}{10q},
\end{eqnarray*}
which
 proves the last inequality in (\ref{2.5h}).

The first inequality in (\ref{2.5h}) follows from
(\ref{3.2}) and (\ref{3.1pi}):
\begin{eqnarray*}
\frac{|\xi^j_{s,k}-\xi^j_{s,k-1}|}{d(\xi^j_{s,k},K^j)}
&\ge&
\frac{|\psi_j(\xi^j_{s,k})-\psi_j(\xi^j_{s,k-1}|)}{16(e^{s/\om_j}-1)}\\
&\ge&
\frac{1}{16}\frac{2}{\pi}\frac{2\pi}{m_j}\frac{\om_j m}{c_2 c q}\ge\frac{1}{c^2q}.
\end{eqnarray*}

For $\z\in D^j_{s,k}$ and $\z^*:=\psi_j((e^{s/\om_j}-1)\phi_j(\z)/|\phi_j(\z)|)$
by virtue of (\ref{3.2}), we have
$$
d(\z,I^j_{s,k})\le |\z-\z^*|<\frac{1}{4} d(\z^*,K^j).
$$
Therefore, by (\ref{2.5h}),
\begin{eqnarray*}
|\z-\xi^j_{s,k}|&\le&|\z-\z^*|+|\z^*-\xi^j_{s,k}|\\
&\le& \frac{1}{4}\left(|\z^*-\xi^j_{s,k}|+ d(\xi^j_{s,k},K^j)\right)
+|\z^*-\xi^j_{s,k}|\\
&\le& \left(\frac{5}{40}+\frac{1}{4}\right)d(\xi^j_{s,k},K^j),
\end{eqnarray*}
which yields (\ref{2.2du}).

Moreover, (\ref{2.1pu}) implies for $\tau=\phi(\z)$ and
$\z\in D^j_{s,k}$,
$$
|\psi_j'(\tau)|\ge\frac{1}{4}\frac{d(\z,K^j)}{e^{s/\om_j}-1}
\ge\frac{1}{8}\frac{d(\xi^j_{s,k},K^j)}{e^{s/\om_j}-1}
$$
which yields (\ref{2.1du}) as follows
\begin{eqnarray*}
\sigma(D^j_{s,k})&=&
\int\int_{\tilde{D}^j_{s,k}}|\psi_j'(\tau)|^2d\sigma(\tau)
\ge\frac{1}{64}\frac{d(\xi^j_{s,k},K^j)^2}{(e^{s/\om_j}-1)^2}
\sigma(\tilde{D}^j_{s,k})\\
&=&
\frac{1}{64}\frac{d(\xi^j_{s,k},K^j)^2}{e^{s/\om_j}-1}
|\tilde{I}^j_{s,k}| \ge
\frac{1}{64}\frac{2\pi}{m_j}\frac{\om_j}{c_2s}d(\xi^j_{s,k},K^j)^2
\ge\frac{d(\xi^j_{s,k},K^j)^2}{c^2}.
\end{eqnarray*}

 \hfill$\Box$

Since by Lemma \ref{lem2.1} and (\ref{2.4pu}) for $z\in K_{10s}$,
$$
\frac{2qd^j_{s,k}}{d(z,I^j_{s,k})}\le\frac{2}{5},
$$
according to (\ref{2.11}) we have
\begin{eqnarray*}
&&|\log|P_n(z)|-ng_\Om(z)-n\log \kap(K)|\\
&\ole& q^2+q^3\left(\frac{2}{5}\right)^{q}
\sum_{j=1}^\nu\sum_{k=1}^{m_j}\left(\frac{d^j_{s,k}}{d(z,I^j_{s,k})}
\right)^2.
\end{eqnarray*}
Furthermore, by Lemma \ref{lem2.1}, (\ref{2.1du}) and (\ref{2.4pu}), 
for $z\in K_{10s}$,
$$
\sum_{j=1}^\nu\sum_{k=1}^{m_j}\left(\frac{d^j_{s,k}}{d(z,I^j_{s,k})}
\right)^2\ole
\sum_{j=1}^\nu\sum_{k=1}^{m_j}\frac{\sigma(D^j_{s,k})}{d(z,D^j_{s,k})^2}
\ole\int\int_{D_s}\frac{d\sigma(\z)}{|\z-z|^2}.
$$
Since by the Loewner inequality (see \cite[p. 27, Lemma 2.5]{andbla}),
$
d(z,K_s)\oge s^2 ,
$
using polar coordinates with center at $z$, we obtain
$$
\int\int_{D_s}\frac{d\sigma(\z)}{|\z-z|^2}\ole\log\frac{\mb{diam}(K_{s_0})}{d(z,K_s)},
\ole\log\frac{s_0}{s},
$$
which yields
\beq\label{2.8du}
|\log|P_n(z)|-ng_\Om(z)-n\log \kap(K)|\ole q^2+2^{-q}\log m.
\eeq

{\bf Proof of Theorem \ref{th1}}.
Without loss of generality, we assume that $n$ is sufficiently large.
First, let $n=mq$, where $q=q_m:=\lfloor 2\log\log m\rfloor.$ Then
by (\ref{2.8du}), for $z\in K_{10cq/m}$,
$$
|\log|P_n(z)|-ng_\Om(z)-n\log\kap(K)|\ole(\log\log m)^2
\asymp (\log\log n)^2.
$$
Thus, the maximum principle implies that for $P_n^*:= P_n\kap(K)^n,
s=cq/m=cq^2/n$ and $z\in\C\setminus U_{10s}$,
$$
\exp\left( ng_\Om(z)-c_3(\log\log m)^2\right)\le|P_n^*(z)|
\le \exp\left( ng_\Om(z)+c_3(\log\log m)^2\right).
$$
Therefore,
$$
||P_n^*||_K\le ||P_n^*||_{K_{10s}}\le\exp\left(c_3(\log\log n)^2\right).
$$
At the same time, for $z\in\C\setminus U_\de,\de=3c_3(\log\log n)/n$,
$$
|P_n^*(z)|\ge \exp\left(2c_3(\log\log m)^2\right)>||P_n^*||_K,
$$
which proves (\ref{1.1r}) for $n=mq_m$.

For arbitrary (sufficiently large) $n$ we find $m\in\N$ such that
$$
mq_m\le n<(m+1)q_{m+1}.
$$
Since
$$
1\le\frac{n}{mg_m}\le\frac{(m+1)q_{m+1}}{mq_{m}}\to1\quad\mb{as }n\to\infty,
$$
we obtain
$$
s_n(K)\le s_{mq_m}(K)\ole\frac{(\log\log m)^2}{mq_m}\ole\frac{(\log\log n)^2}{n},
$$
which completes the proof of (\ref{1.1r}).

\hfill$\Box$

From now on we assume that each $K^j,j=1,\ldots,\nu$ is a quasidisk.
Since $\partial\Om^j$ consists of two quasiconformal curves, $\phi_j$
can be extended to a $Q_j$-quasiconformal homeomorphism of a neighborhood
of $\ov{\Om^j}$ to a neighborhood of $\ov{A^j}$ with some $Q_j\ge 1$.
Therefore, repeating the proof of \cite[p. 29, Theorem 2.7]{andbla}
 we can establish the following result.
\begin{lem}\label{lem3.1}
Let $\z_k\in\Om^j, w_k:=\phi_j(\z_k), k=1,2,3$. Then:

(i) the conditions $|\z_1-\z_2|\le C_1|\z_1-\z_3|$ and
  $|w_1-w_2|\le C_2|w_1-w_3|$ are equivalent; besides, the constants 
	$C_1$ and $C_2$ are mutually dependent and depend on $Q_j$ and $K$;
	
	(ii) 
if $|\z_1-\z_2|\le C_1|\z_1-\z_3|$, then
$$
\frac{1}{C_3}\left|\frac{w_1-w_3}{w_1-w_2}\right|^{1/Q}\le
\left|\frac{\z_1-\z_3}{\z_1-\z_2}\right|\le
C_3\left|\frac{w_1-w_3}{w_1-w_2}\right|^{Q},
$$
where $Q:=\max_{j}Q_j$ and $C_3=C_3(C_1,Q,K)>1$.
\end{lem}
For $\z\in \Om^j$ let $\tilde{\z}:=\psi_j(\phi_j(\z)/|\phi_j(\z)|)$.
We fix $\z^*_j\in K^j_{s_0}$. By Lemma \ref{lem3.1}, for
$\z\in K_s^j$ and $z\in K^j_{10s}, s< s^*$,
\beq\label{3.11}
\frac{d(\z,K^j)}{|\z-z|}\le\left|\frac{\z-\tilde{\z}}{\z-z}\right|
\ole \left(\frac{s}{|\phi_j(\z)-\phi_j(z)|}\right)^{1/Q},
\eeq
\beq\label{3.12}
d(\z,K^j)\le|\z-\tilde{\z}|\asymp
\left|\frac{\z-\tilde{\z}}{\z-\z^*_j}\right|
\ole\left|\frac{\phi_j(\z)-\phi_j(\tilde{\z})}{\phi_j(\z)-\phi_j(\z^*_j)}\right|^{1/Q}
\asymp s^{1/Q}.
\eeq
Furthermore, by virtue of Lemma \ref{lem2.1}, for $z\in K_{10s}$ and $q\in\N$,
\begin{eqnarray*}
&&\int_{K_s^j}\frac{d(\z,K^j)^q}{|\z-z|^{q+1}}|d\z|=
\sum_{k=1}^{m_j}\int_{I_{s,k}^j}\frac{d(\z,K^j)^q}{|\z-z|^{q+1}}|d\z|\\
&\oge&\sum_{k=1}^{m_j}\frac{d(I_{s,k}^j,K^j)^q|I_{s,k}^j|}{d(z,I_{s,k}^j)^{q+1}}
\oge \sum_{k=1}^{m_j}\left(\frac{d_{s,k}^j}{d(z,I_{s,k}^j)}\right)^{q+1}.
\end{eqnarray*}
Therefore, according to (\ref{2.11})
for $z\in K^{r}_{10s}, r=1,\ldots,\nu$,
\begin{eqnarray*}
&&|\log|P_n(z)|-ng_\Om(z)-n\log\kap(K)|\\
&\ole&1+\sum_{j=1,j\neq r}^\nu\int_{K^j_s}d(\z,K^j)^q|d\z|
+
\int_{K_s^r}
\frac{d(\z,K^r)^q}{|\z-z|^{q+1}}|d\z|.
\end{eqnarray*}
Let $q:=\lfloor 2Q \rfloor$ and $ w=\phi_r(z)$. Then 
by (\ref{2.1pu}), (\ref{3.11}) and (\ref{3.12}),
\begin{eqnarray*}
\int_{K^j_s}d(\z,K^j)^q|d\z|&\ole&\frac{1}{s}
\int_{|\tau|=e^{s/\om_j}}d(\psi_j(\tau),K^j)^{q+1}|d\tau|\\
&\ole& s^{(q+1)/Q-1}\ole 1.
\end{eqnarray*}
as well as
\begin{eqnarray*}
\int_{K_s^r}
\frac{d(\z,K^r)^q}{|\z-z|^{q+1}}|d\z|&\ole&
\frac{1}{s}\int_{|\tau|=e^{s/\om_r}}\left(\frac{d(\psi_r(\tau),K^r)}
{|\psi_r(\tau)-\psi_r(w)|}\right)^{q+1}|d\tau|\\
&\ole&
\frac{1}{s}\int_{|\tau|=e^{s/\om_r}}\frac{s^{(q+1)/Q}}{|\tau-w|^{(q+1)/Q}}
|d\tau|\ole 1.
\end{eqnarray*}
Thus,
the maximum principle implies that for $P^*_n(z):=P_n(z)\kap(K)^n$
and $s=cq/m=cq^2/n, n=qm,m>m_0$,
\beq\label{3.20}
\frac{1}{c_4}e^{ng_\Om(z)}\le|P_n^*(z)|\le c_4e^{ng_\Om(z)},\quad
z\in\C\setminus U_{10s}.
\eeq

{\bf Proof of Theorem \ref{th2}}.
Without loss of generality, we assume that $n$ is sufficiently large.
First, let $n=mq$ be such that (\ref{3.20})
holds. We have
$$
||P^*_n||_K\le ||P^*_n||_{K_{10s}}\le c_4\exp\left(10cq^2\right).
$$
If we let
$$
c_5:=2\log c_4+11cq^2,
$$
then for $z\in\C\setminus U_{c_5/n}$, we get
$$
|P_n^*(z)|\ge\frac{1}{c_4}\exp(c_5)>c_4\exp\left(10cq^2\right),
$$
which shows that $s_n(K)\le c_5/n$.

If $mq<n\le m(q+1)$, then
$$
s_n(K)\le\frac{c_5}{mq}\le \frac{2c_5}{m(q+1)}\le\frac{2c_5}{n}.
$$
Hence, in both cases we have (\ref{1.4}).

 \hfill$\Box$

{\bf Acknowledgements.}
The author is  grateful to   M. Nesterenko
 for his helpful comments.

V. V. Andrievskii

 Department of Mathematical Sciences

 Kent State University

 Kent, OH 44242

 USA

e-mail: andriyev@math.kent.edu


\begin{thebibliography}{99}


\bibitem{ahl} Ahlfors L. V. (1966)

Lectures on Quasiconformal
Mappings,  Van Nostrand, Princeton, N. J.


\bibitem{ahl1}
 Ahlfors L. V. (1979)

Complex Analysis,  McGraw-Hill, Inc., New York.
	
\bibitem{and00}
 Andrievskii V. V. (2000)

 On the approximation of a continuum by lemniscates,
J. Approx. Theory 105, 292-304.


 \bibitem{and151}
 Andrievskii V. V. (2016)

 Chebyshev polynomials on a system of continua,
 Constr. Approx.  43, 217-229.
	

 \bibitem{andbla}  Andrievskii V. V. and  H. -P. Blatt (2002)

Discrepancy of Signed Measures and
      Polynomial Approximation,
       Springer-Verlag, Berlin/New York.
			
	\bibitem{andnaz}
	 Andrievskii V. V. and   F. Nazarov (2018)
	
	On the Totik-Widom property for a quasidisk,
	arXiv:1802.06948.
	
				
	\bibitem{gai}
			Gaier D. (1999)
			
			The Faber operator and its boundedness,
			J. Approx. Theory 101, 265-277.
			

			\bibitem{geh}
 Gehring F. W.  (1982)

Characteristic Properties of  Quasidisks,
Les Presses De L'Universit\'e De Montr\'eal, Montr\'eal.	

\bibitem{kit}
 Laiyi Z. (2004)

On the degree of convergence of lemniscates in finite connected domains,
J. Approx. Theory 131, 185-195.

\bibitem{lehvir}
Lehto O. and K. I. Virtanen (1973)

 Quasiconformal Mappings in the
Plane, 2nd ed., Springer-Verlag, New York.


\bibitem{pri}
Pritsker I.  (2011)

Equidistribution of points via energy,
Ark. Mat. 49, 149-173.

\bibitem{ran}
  Ransford T. (1995)
		
		Potential Theory in the Complex Plane,
   Cambridge University Press, Cambridge.


\bibitem{saftot}
 Saff E. B. and    V. Totik  (1997)

    Logarithmic Potentials with External
Fields, Springer-Verlag,  New York/Berlin.


\bibitem{sic}
Siciak J. (1967)

Degree of convergence of some sequences in the conformal mapping
theory, Coll. Math. 16, 49-59.


\bibitem{tot12}
Totik V. (2012)

 Chebyshev polynomials on a system of curves,
Journal D'Analyse Math\'ematique 118, 317-338.

\bibitem{tot13}
Totik V. (2014)
 
 Chebyshev polynomials on compact sets,
Potential Anal. 40, 511-524.




\bibitem{wid}
 Widom H. (1969)

 Extremal polynomials associated with a
system of curves in the complex plane,
 Adv. in Math.  3, 127--232.

\end{thebibliography}
\end{document}